\documentclass[conference]{IEEEtran}
\IEEEoverridecommandlockouts
\usepackage[keeplastbox]{flushend}
\usepackage{cite}
\usepackage{amssymb,amsfonts}
\usepackage[noend]{algpseudocode}
\usepackage{multicol}
\usepackage{graphicx}
\usepackage{textcomp}
\usepackage{xcolor}
\usepackage{flushend}
\usepackage{amsmath}
\usepackage{fancyhdr}
\usepackage{float}

\usepackage[obeyspaces]{url}
\usepackage[hardware,software]{mymacros}
\usepackage{amsthm}
\usepackage{pgfplots}

\usepackage[vlined,longend,linesnumbered,ruled]{algorithm2e}
\providecommand{\dontprintsemicolon}{\DontPrintSemicolon}
\providecommand{\DontPrintSemicolon}{\dontprintsemicolon}

\usepackage[list=true]{caption,subcaption} 
\usepackage{pgfplots}
\definecolor{mpiblue}{HTML}{006AA9}
\definecolor{mpigreen}{HTML}{5C871D}
\pgfplotscreateplotcyclelist{mylist1}{%
{black,solid},
{blue,dashed},
{mpigreen,dash pattern = on 3pt off 1pt on 1pt off 1pt},
{red,dash pattern = on 2pt off 2pt on 1pt off 2pt on 1pt off 2pt},
{cyan,densely dashed},
{mpiblue,densely dotted},
{magenta,dotted},
}
\usepackage{pgfplots}
\definecolor{mpiblue}{HTML}{006AA9}
\definecolor{mpigreen}{HTML}{5C871D}
\pgfplotscreateplotcyclelist{mylist1}{%
{black,solid},
{blue,dashed},
{mpigreen,dash pattern = on 3pt off 1pt on 1pt off 1pt},
{red,dash pattern = on 2pt off 2pt on 1pt off 2pt on 1pt off 2pt},
{cyan,densely dashed},
{mpiblue,densely dotted},
{magenta,dotted},
}
\pgfplotsset{every axis plot/.append style={line width=.75pt}}

\usepackage{eso-pic}

\def\BibTeX{{\rm B\kern-.05em{\sc i\kern-.025em b}\kern-.08em
    T\kern-.1667em\lower.7ex\hbox{E}\kern-.125emX}}
    
\fancypagestyle{specialfooter}{%
  \fancyhf{}
  
  \fancyfoot[C]{}
}
  \usepackage{bm}
\usepackage{tikz}
\usetikzlibrary{shapes.geometric, arrows}
\tikzstyle{startstop} = [rectangle, rounded corners, minimum width=3cm, minimum height=1cm,text centered, draw=black, fill=white!30]
\tikzstyle{io} = [trapezium, trapezium left angle=70, trapezium right angle=110, minimum width=3cm, minimum height=1cm, text centered, draw=black, fill=white!30]
\tikzstyle{process} = [rectangle, minimum width=3cm, minimum height=1cm, text centered, draw=black, fill=white!30]
\tikzstyle{decision} = [diamond, minimum width=3cm, minimum height=1cm, text centered, draw=black, fill=white!30]
\tikzstyle{arrow} = [thick,->,>=stealth]

  \def\@fnsymbol#1{\ensuremath{\ifcase#1\or *\or \dagger\or \ddagger\or
   \mathsection\or \mathparagraph\or \|\or **\or \dagger\dagger
   \or \ddagger\ddagger \else\@ctrerr\fi}}
 \newcommand{\ssymbol}[1]{^{\@fnsymbol{#1}}}   
 \makeatother

\begin{document}
\title{Data Assimilation: Two Different Perspectives Based on the Initial-Condition Dependence}

\author{\textbf{Mohammad N. Murshed}\(\ssymbol{4}\)\ , \textbf{Zarin Subah}\(\ssymbol{1}\)\ , and \textbf{M. Monir Uddin}\(\ssymbol{4}\)\
\\ \\ \(\ssymbol{4}\)\ Department of Mathematics and Physics,\\ North South University, Dhaka, Bangladesh.\\
\(\ssymbol{1}\)\ Institute of Water and Flood Management,\\ Bangladesh University of Engineering and Technology, Dhaka, Bangladesh.}


\IEEEoverridecommandlockouts
\maketitle
\thispagestyle{specialfooter}

\IEEEpubidadjcol
%
\begin{abstract}
Data Assimilation (DA) is a computational tool that uses value from the model and the real measurement to arrive to an optimally acceptable value. Rather, this technique relies on the idea of Kalman gain. We point out that DA has two different perspectives based on the type of problem. In this paper, we look into two problem types: one that does not rely on the initial condition, and the other that is initial condition dependent. Data Assimilation is demonstrated on two examples: runoff monitoring and forecasting in the city of Dhaka (initial condition independent) and convection in the atmosphere (initial condition dependent). We show that standard DA works well for problems with no initial condition dependence and piecewise DA is to be utilized when the problem has initial condition dependence. In the first example, we exploited standard Data Assimilation to arrive at values that are more realistic than the ones from the model and the observations. The second example is where we devised a method to find the dynamics of the system in a piecewise manner in Data Assimilation framework and noticed that the data assimilated dynamics (even due to noisy initial condition) is in good agreement with the true dynamics for a reasonable extent of time in future.  \\
\end{abstract}

%
%
\section{Introduction}
\begin{figure}
\centering
\begin{tikzpicture}[node distance=2cm]
\node (start)[startstop]{Problem Type};
\node (pro1a) [process, below of=start, yshift= 0.5 cm, xshift=-2.3 cm] {Initial condition independent};
\node (pro2a) [process, below of=start, yshift= 0.5 cm, xshift= 2.3 cm] {Initial condition dependent};
\node (out1) [process, below of=pro1a,yshift= 0.4 cm, xshift = 0  cm ] {Standard DA};
\node (out2) [process, below of=pro2a, yshift= 0.4 cm, xshift= 0 cm] {Piecewise DA};
\draw [arrow] (start) -- (pro1a);
\draw [arrow] (start) -- (pro2a);
\draw [arrow] (pro1a) --(out1);
\draw [arrow] (pro2a) --(out2);
\end{tikzpicture}\\
\caption{DA Perspective based on the initial condition dependence of the problem}
\label{fig:Arch}
\end{figure}
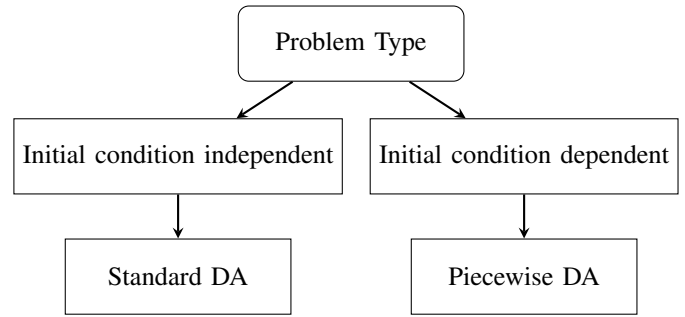

Data Assimilation (DA) is a numerical method of weighing observation and values from the model to give the best possible estimate of a quantity of interest in a problem,\cite{evensen2009data, kutz2013data,ghil1991data}. The weighing is usually done based on the uncertainties associated with the observation and the model. DA is sometimes thought to be related to the notions from Bayes' theorem. It enjoys popularity in the field of weather prediction where it can be used to reconstruct and somewhat predict the behavior of the atmosphere over time. Moreover, this technique is capable of extrapolating the dynamics of atmosphere-ocean interaction beyond a particular region for which the data is available, \cite{derber1989global,robinson2001data}. \\ \\  
We consider a problem where we seek to use the value of a quantity as an input to compute the value of another quantity (the output). For instance, \(v\) can be a function of \(u\),
$$ v=f(u) $$ hence, the value of \(u\) must be set beforehand. The value of \(u\) is accessible from either a model or from real measurements. Although, it is easy to produce values from a model, using values from a model makes it less real. On the other hand, observations are real, but, not always available (sparse) for a given problem. This makes it necessary to take into account both the assumptions in a model and the physical aspects of the problem to estimate the current and future state of the variables of the system accurately, \cite{bouttier2002data,nichols2003data}. This process of melding dynamics and data is a robust methodology which helps to find efficient and realistic estimation.\\ \\
In this work, we show how to apply and interpret DA as per the dependence of the problem on the initial condition, depicted in Figure \ref{fig:Arch}. Standard DA works fine for problems with no initial condition dependence, but DA needs to be applied in a piecewise manner when the problem is initial condition dependent. This paper is organized as follows: Section \ref{BACK} is a review of the basics of the standard DA. The DA perspective to deal with initial condition dependent problem is elaborated in Section \ref{NPKG}. The results are shown in Section \ref{NR} and a summary provided in Section \ref{CFW}. 

\section{Background}
\label{BACK}
Assume that we have two values for today's temperature: one is experimental, \(T_e\), and another from a model, \(T_m\). If the uncertainties from the model, \(\sigma_{m}\), and the experiment, \(\sigma_{e}\), are provided, DA offers a way to consider data from both the model and experiment to reach an optimal value, \(T_{DA}\), for the temperature,
\begin{equation}
 T_{DA} = \beta_{m} T_{m} + \beta_{e}T_{e},
 \end{equation}  
where \(\beta_{m}\) and \(\beta_{e}\) are weighing factors used for the value from the model and the one from the experiment, respectively. The sum of these factors equals to 1. It is worth noting that the value from the model and the experimental value are unbiased which implies that \(\sigma_{m} \sigma_{e} \approx 0\). \\ \\
The objective of DA is to look for \(\beta_{m}\) and \(\beta_{e}\) to optimize the uncertainty in \(T_{DA}\), \(\sigma_{DA}\). Optimization, \cite{ghil1991data}, results in the following optimal weights and the uncertainty for the data-assimilated value:
\begin{equation}
\beta_{m} = \frac{\sigma_{e}^2}{\sigma_{m}^2+ \sigma_{e}^2} 
\end{equation}
\begin{equation}
\beta_{e} = \frac{\sigma_{m}^2}{\sigma_{m}^2+ \sigma_{e}^2} 
\end{equation}
\begin{equation}
\sigma_{DA}^{2} = \frac{\sigma_{m}^2\sigma_{e}^2}{\sigma_{m}^2+ \sigma_{e}^2} < \sigma_{m}^{2}, \sigma_{e}^{2}.
\end{equation}
The uncertainty in the data-assimilated value would be less than that of the value from the model and that of the one from the experiment. Such formulations yield an optimal value of the temperature that is based on the output from the model and the data from experiment. The additional information required are the errors in the model and the experiment. 

\section{Perspective on Data-Assimilation based on the initial condition of the problem}
\label{NPKG}
Here, we discuss the role of Kalman gain in Data Assimilation and then show how DA can be employed in a setting where initial condition has influence on the problem.\\ \\
Kalman gain leverages the expression of the data-assimilated value in terms of the model value and the experimental value by,
\begin{equation}
\label{KG}
T_{DA,k+1}=T_{m,k+1} + W_{k+1}(T_{e,k+1} -T_{m,k+1})
\end{equation}
where 
$$ W_{k+1} = \frac{\sigma_{m,k+1}^2}{\sigma_{m,k+1}^2+\sigma_{e,k+1}^2}. $$
\(W\) is the Kalman gain, \(k\) indicates time and \(\sigma^2\) means the variance. Note that if \(\sigma_e^{2}\) goes to zero, then \(W\) = 1 which implies that \(T_{DA}=T_e\).\\ \\
For an initial condition dependent problem, there may be a small amount of noise in the initial condition, but, such noise is not supposed to give a spurious solution. In case the dynamics is not accurate just because of the some noise in the initial condition, data-assimilation enables increasing the accuracy of the dynamics upto a certain time instant in future. We propose a pseudocode that can predict the solution reliably for a much longer time: \\ \\
1. Find the dynamics of the quantity of interest over a time span of \(N_t\) temporal nodes: \(T_{m,1},T_{m,2},T_{m,3},........, T_{m,N_t}\). (This can be done easily via \textbf{ode45} which is a function on MATLAB to solve differential equation of the given problem). \\ \\
2. Record real measurements every \(q\) time steps: \(T_{e,q+1},T_{e,2q+1},T_{e,3q+1}\). \\ \\
3. Use \textbf{ode45} to find the dynamics from the initial time, \(t_{1}\), upto the time of first measurement, \(t_{q+1}\). \\ \\
4. Assimilate \(T_{m,q+1}\) and \(T_{e,q+1}\) by Eq (\ref{KG}). \\ \\
5. This data-assimilated value would then be used as the initial condition for next round of simulation from \(t=q+2\) to \(t=2q+1\). \\ \\
6. Step 5 can be repeated till we get to the final time. \\ 

\section{Numerical Results}
\label{NR}
We apply the standard DA on 'Runoff in Dhaka city'  (independent of the initial condition) and then illustrate the utility of piecewise DMD on 'Convection in the atmosphere' (dependent on initial condition).  Note that all the results are all produced by using MATLAB (R2020a on an Intel CORE i5 processor with 8 GB 1600 MHz DDR3 memory).
\begin{figure}
\centering
\includegraphics[scale=0.37]{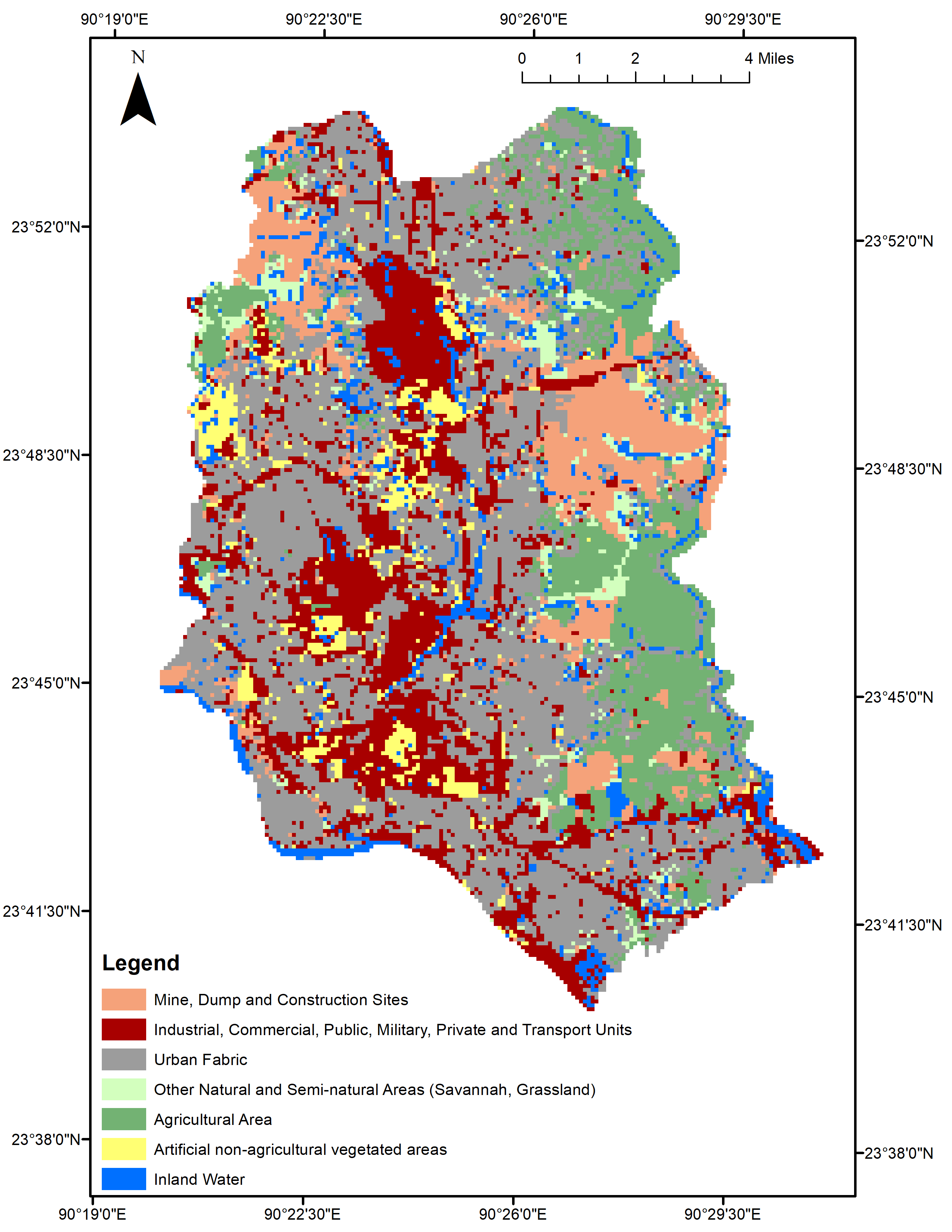}
\caption{Land Use Type of Dhaka City}
\label{g1}
\end{figure}

\subsection*{\textbf{Runoff monitoring and forecasting in Dhaka city}}
\begin{figure}
\centering
\includegraphics[scale=0.53]{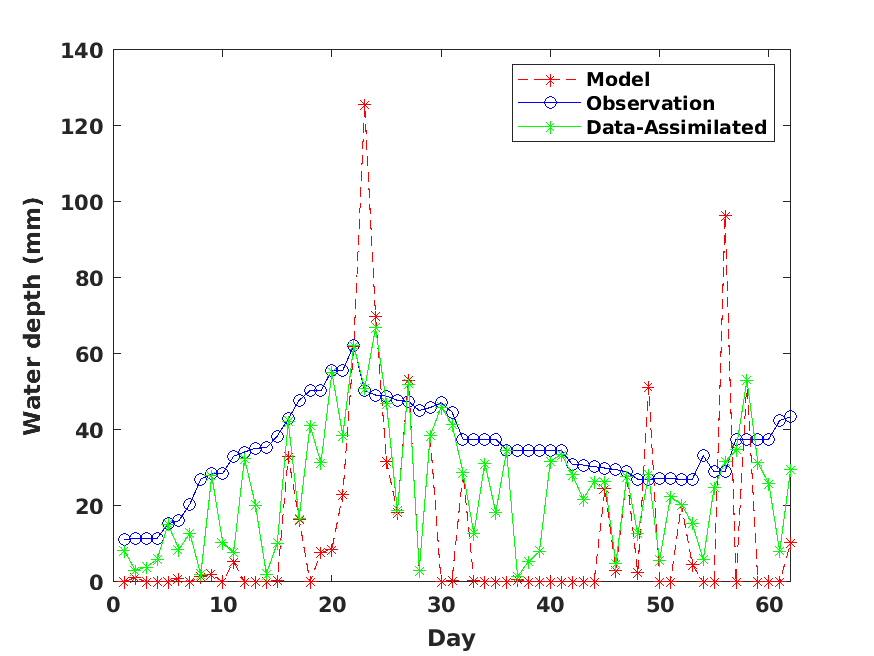}
\caption{Runoff values from the model, measurements and data-assimilation}
\label{rp1}
\end{figure}

\begin{figure}
\centering
\includegraphics[scale=0.53]{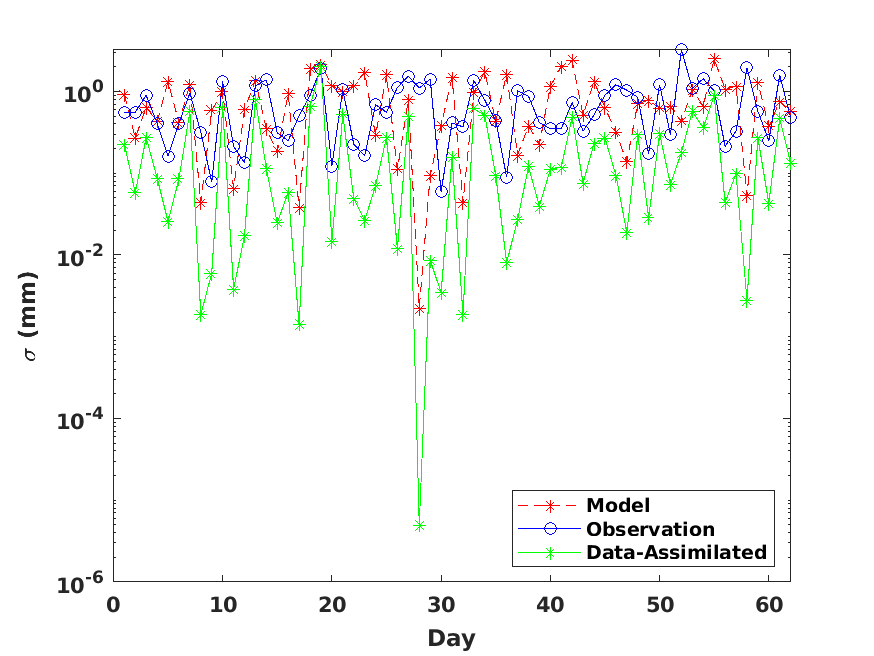}
\caption{Variation in uncertainty present in model, measurement and data-assimilation over time}
\label{rp2}
\end{figure}

Urban flooding in Dhaka city creates inevitable problems for the city dwellers due to short duration heavy rainfall events, poor drainage system management and encroachment of the wetlands and waterbodies, \cite{braun2011floods}. Though severe urban flooding events have occurred in recent decades, there are limited rainfall gauges and water level monitoring systems available in the city which create challenges for providing accurate runoff prediction. For this reason, physical models can attribute anomalies in predicting short term flooding. In this scenario, data assimilation can be a useful technique for providing a more reliable prediction and updating the model state, \cite{fava2018approach}. We focused on the flood that occurred in 2007 in order to understand the usefulness of DA in the monitoring and forecasting runoff of Dhaka city. \\ \\
Rainfall-runoff estimation for this study is done using Soil Conservation Services and Curve Number (SCS-CN)  and Geographic Information System (GIS) approach. Integration of SCS-CN model into remote sensing and GIS systems extends the applicability of the model to complex variability with high spatiotemporal variability in land use pattern. In this study, processing of the land use and land cover data of Dhaka city in 2007 (obtained from Earth Observation for Sustainable Development by  European Space Agency) by GIS gives seven types of land use patterns, Figure \ref{g1}. From the table given by Hong and Adler (2008), modified from USDA (1986) handbook, CN values for Dhaka city are extracted from land use patterns and Hydrologic soil group data for Dhaka (collected from NASA Earthdata), [11], [12]. \\ \\
SCS-CN model estimates runoff (\(Q\)) as,
\begin{equation}
Q=\frac{(P-I_{a})^{2}}{P-I_{a} +S}
\end{equation}
where \(P\) is the rainfall depth, \(I_a\) initial abstraction, \(S\) potential maximum retention after runoff begins, \cite{cronshey1986urban}. \(S\) \cite{caviedes2012influence} is calculated as per,
\begin{equation}
S = \frac{25400}{CN} - 254.
\end{equation} 
Weighted curve number is computed using area specific curve numbers for different land use and soil condition. It is important to note that curve number values, in this particular context, are also influenced by Antecedent Moisture Condition that relates to the rainfall values (July to August, 2007) obtained from Bangladesh Water Development Board. This is how we obtained the runoff values from the model and the observed data are collected from Dhaka Water Supply and Sewerage Authority.\\ \\
Both the model values (red) and the observations (blue) are displayed in Figure \ref{rp1}. It is imperative to note that we introduced the uncertainty in the model state and the observations by considering absolute value of normally distributed random numbers. Application of DA on the run-off data from the model and the observed data results in a more realistic, reliable value for each day, also shown in Figure \ref{rp1}. The uncertainty in data-assimilated values appears to be much lower than both the uncertainty in the values from the model and the uncertainty in the observations, Figure \ref{rp2}.

\subsection*{\textbf{Convection in the Atmosphere}}
\begin{figure}
\centering
\includegraphics[scale=0.4]{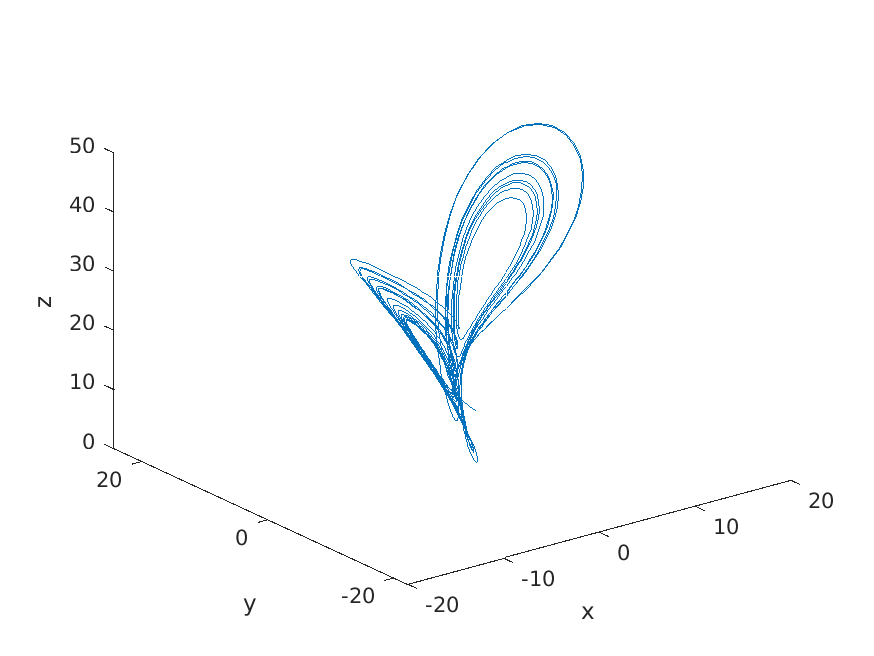}
\caption{Butterfly effect in convection in the atmosphere}
\label{lm1}
\end{figure}

\begin{figure}
\centering
\includegraphics[scale=0.5]{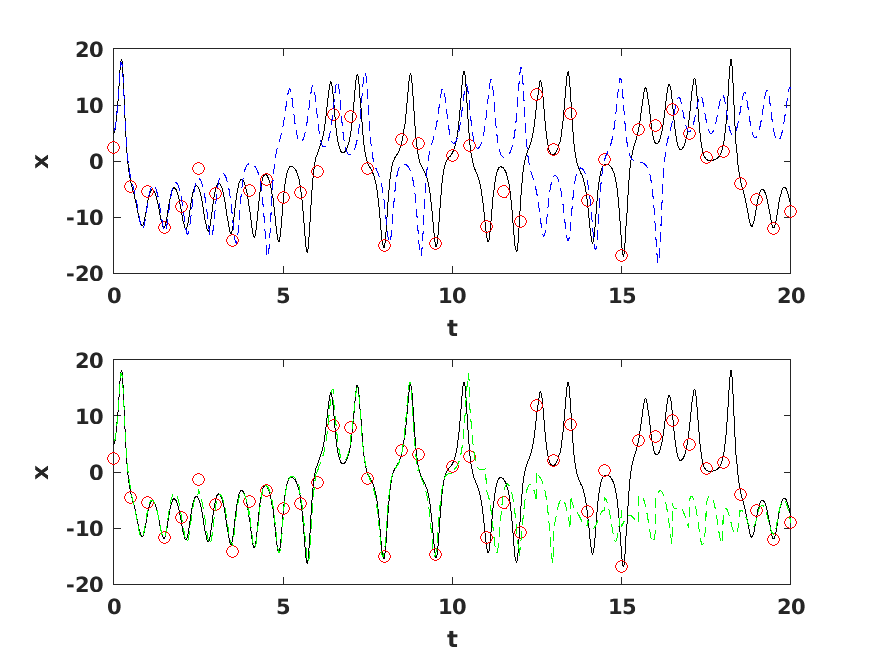}
\caption{Top panel: true dynamics (black line),  measurements (circles) and dynamics from noisy initial condition (blue dotted line); Bottom panel: true dynamics (black line),  measurements (circles) and data-assimilated solution (green dotted line)}
\label{lm2}
\end{figure}

\begin{figure}
\centering
\includegraphics[scale=0.6]{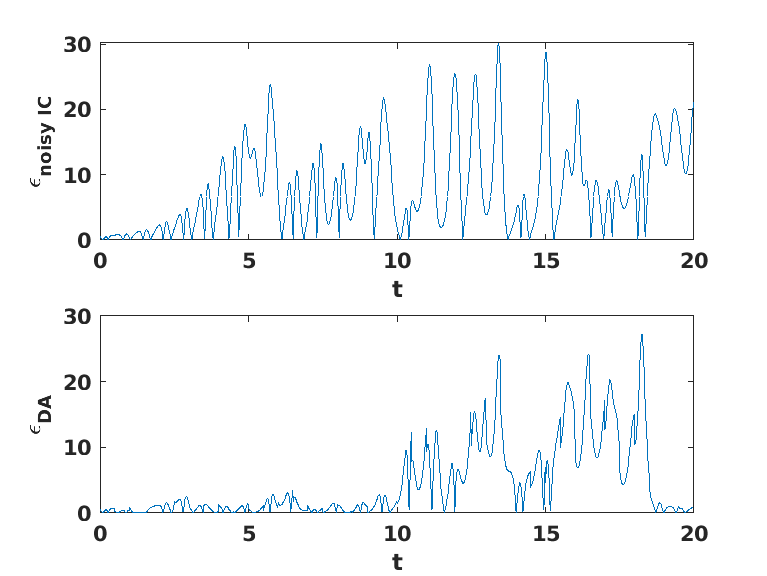}
\caption{Top panel: error in the dynamics from noisy initial condition; Bottom panel: error in the data-assimilated solution}
\label{lm3}
\end{figure}

The convection in the atmosphere is governed by the Lorenz model which is a dynamical system that captures the random behavior of a particle, Figure \ref{lm1}. This phenomenon is also called the 'butterfly' effect. The model for this problem is,
\begin{equation}
\frac{dx}{dt}=a_{1}(y-x)
\end{equation}
\begin{equation}
\frac{dy}{dt}=a_{2}x - y - xz
\end{equation}
\begin{equation}
\frac{dz}{dt}=xy- a_{3}z
\end{equation}
where \(a_{1}=10,a_{2}=28\) and \(a_{3}=8/3\). With these parameters, this dynamical system is excited by the initial condition \((x_{0},y_{0},z_{0})= (5,5,5)\) via \textbf{ode45}. We note the dynamics of \(x\) with time for 20 s. This is how we get the values from the model and we record measurements every 0.5 s.\\ \\
In Figure \ref{lm2}, the black line represents the true solution of this system and the red circles are the measurements taken every 0.5 s. The blue dotted line, in the top panel, is the dynamics of the system when noise is added to the initial condition whereas the green dotted line, in the bottom panel, shows the data-assimilated solution. The noise in the initial condition generates a solution that deviates significantly from the true solution. On the other hand, the data-assimilated solution agrees well with the true solution for a much longer time (upto \(t=12.5 \ s\)). An error analysis based on the metrics,
\begin{equation}
\epsilon_{noisy\  IC}= | x_{true} - x_{noisy\  IC}| 
\end{equation}
\begin{equation}
\epsilon_{DA}= | x_{true} - x_{DA}|, 
\end{equation}
reveals that the solution from the noisy initial condition sees growing error from the very beginning, Figure \ref{lm3} (top panel), and the uncertainty in the data-assimilated solution stays low for an appreciable extent of time and then increases a little bit over time, Figure \ref{lm3} (bottom panel). This proves that data-assimilation can preserve the true structure of the solution for a reasonable time window.

\section{Conclusion and Future Work}
\label{CFW}
In this work, we apply DA on two different test cases: runoff monitoring and forecasting in Dhaka city and convection in atmosphere. In the former case, conventional DA approach is enacted and the resulting data-assimilated values appear to be much promising compared to the ones from the model and measurements. In the latter, we show the use of piecewise DA to handle problems having noisy initial condition. Thus, DA can be put into use in two different ways based on the type of the problem. In future, we plan to use three-/four-dimensional variational analysis method for a better understanding of the oceanic-atmospheric dynamics in the Earth system model.



\bibliography{Turb_p1_ref}
\bibliographystyle{IEEEtran}      

\end{document}